\documentclass{article}

\usepackage{amsmath}  
\usepackage{amssymb}  
\usepackage{dsfont}
\usepackage{amsthm}   
\usepackage{amsfonts}
\usepackage{mathrsfs}

\usepackage{mathrsfs}  
\usepackage{paralist}                       

\swapnumbers
\newtheorem{theorem}{Theorem}[section]
\newtheorem{lemma}[theorem]{Lemma}
\newtheorem{proposition}[theorem]{Proposition}

\newtheorem{coro}[theorem]{Corollary}

\newtheorem*{exampleso}{Examples}

\newcommand*{\hfillplus}{\hfill\linebreak[3]\hspace*{\fill}}

\author{M.~Hellus, R.~Waldi}
\title{Distribution of weights and a question of Wilf}
\date{\today}

\begin{document}

\maketitle

\begin{abstract}

Let $S$ be a numerical semigroup of embedding dimension $e$ and conductor $c$. The question of Wilf is, if $\#(\mathds N\setminus S)/c\leq e-1/e$.

\noindent In (An asymptotic result concerning a question of Wilf, arXiv:1111.2779v1 [math.CO], 2011, Lemma 3), Zhai has shown an analogous inequality for the distribution of weights $x\cdot\gamma$, $x\in\mathds N^d$, w.\,r. to a positive weight vector $\gamma$:

\noindent Let $B\subseteq\mathds N^d$ be finite and the complement of an $\mathds N^d$-ideal. Denote by $\operatorname{mean}(B\cdot\gamma)$ the average weight of $B$. Then
\[\operatorname{mean}(B\cdot\gamma)/\max(B\cdot\gamma)\leq d/d+1.\]

\begin{itemize}

\item For the family $\Delta_n:=\{x\in\mathds N^d|x\cdot\gamma<n+1\}$ of such sets we are able to show, that $\operatorname{mean}(\Delta_n\cdot\gamma)/\max(\Delta_n\cdot\gamma)$ converges to $d/d+1$, as $n$ goes to infinity.

\item Applying Zhai's Lemma 3 to the Hilbert function of a positively graded Artinian algebra yields a new class of numerical semigroups satisfying Wilf's inequality.

\end{itemize}

\end{abstract}

\section{On the distribution of weights at $\mathds N^d$ with regard to a positive weight vector}

\label{section_1}Let $\gamma_1,\ldots,\gamma_d$ be positive real numbers. The \textbf{weight} of $x=(x_1,\ldots,x_d)\in\mathds R^d$ with regard to the \textbf{weight vector} $\gamma=(\gamma_1,\ldots,\gamma_d)$ is defined as the dot product $x\cdot\gamma=x_1\gamma_1+\ldots+x_d\gamma_d$. We cut $\mathds R^d$ into strips
\[H_i=H_i(\gamma)=:=\{x\in\mathds R^d|i\leq x\cdot\gamma<i+1\}\]
and set
\[h_i:=\#(H_i\cap\mathds N^d), i\in\mathds Z.\]
Then
\[(H_0~\dot{\cup}~H_1~\dot{\cup}~\ldots~\dot{\cup}~H_n)\cap\mathds N^d=\Delta_n=\Delta_n(\gamma):=\{m\in\mathds N^d|m\cdot\gamma<n+1\},\]
hence
\[\#\Delta_n=\sum_{j=0}^nh_j\text{ and }\sum_{m\in\Delta_n}\lfloor m\cdot\gamma\rfloor=\sum_{j=0}^njh_j.\]
Set
\begin{equation}\label{h_ratio}h(n,\gamma):=\frac{\sum_{j=0}^njh_j}{n\sum_{j=0}^nh_j}=\frac{\sum_{m\in\Delta_n}\lfloor m\cdot\gamma\rfloor}{n\cdot\#\Delta_n}\text{ for }n>0.\end{equation}
In the proof of \ref{Proposition} Prop., we will use the following result of Zhai.

For any finite subset $B\subseteq\mathds N^d$ let
\[\operatorname{mean}(B\cdot\gamma):=\frac1{\#B}\sum_{b\in B}b\cdot\gamma\]
the mean weight of its elements w.\,r. to $\gamma$.

\begin{lemma}\label{Zhai} (\cite[Lemma 3]{Z}) Let $B\subseteq\mathds N^d$ be the complement of an Artinian $\mathds N^d$-ideal. Then
\[\frac{\operatorname{mean}(B\cdot\gamma)}{\max(B\cdot\gamma)}\leq\frac d{d+1}.\]
\end{lemma}

Obviously $\mathds N^d\setminus\Delta_n$ is an ideal in $\mathds N^d$. Hence by \ref{Zhai} Lemma
\begin{equation}\label{z_ratio}\operatorname{mean}(\Delta_n\cdot\gamma)\leq\frac d{d+1}\max(\Delta_n\cdot\gamma)<(n+1)\frac d{d+1}.\end{equation}
Since $m\cdot\gamma-1<\lfloor m\cdot\gamma\rfloor\leq m\cdot\gamma$ we get from (\ref{h_ratio}) and (\ref{z_ratio})
\begin{equation}\label{drei}\frac1n\operatorname{mean}(\Delta_n\cdot\gamma)-\frac1n<h(n,\gamma)\leq\frac1n\operatorname{mean}(\Delta_n\cdot\gamma)\leq\frac1n\frac d{d+1}\max(\Delta_n\cdot\gamma).\end{equation}
This implies

\begin{proposition}\label{Proposition} \begin{enumerate} \item[a)] $h(n,\gamma)<\frac{n+1}n\frac d{d+1}$ for $n>0$.

\item[b)] If $\gamma\in\mathds N^d$, then $h(n,\gamma)\leq\frac d{d+1}$ for $n>0$.\end{enumerate}\end{proposition}

\textbf{Proof }\begin{enumerate}\item[a)] is immediate from (\ref{z_ratio}) and (\ref{drei}).

\item[b)] Here $\Delta_n\cdot\gamma\subseteq\mathds N$, hence $\max(\Delta_n\cdot\gamma)\leq n$ and $h(n,\gamma)\leq\frac d{d+1}$ by (\ref{drei}).\end{enumerate}\hfillplus$\square$

However, there are pairs $(n,\gamma)$ such that $h(n,\gamma)>\frac d{d+1}$.

\noindent\textbf{Example} $d=2$, $n=4$, $h(4,(\frac{12}8,\frac{13}8))=\frac{27}{40}>\frac23$.

\vspace{.5cm}

For general $\gamma\in\mathds R_{>0}^d$, we have the following asymptotic result.

\begin{theorem}\label{thm_13}\begin{enumerate}\item[a)] $\lim_{n\to\infty}h(n,\gamma)=\frac d{d+1}$ for each $\gamma\in\mathds R_{>0}^d$.

\item[b)] For $\gamma=(1,\ldots,1)$ even $h(n,\gamma)\equiv\frac d{d+1}$.\end{enumerate}\end{theorem}

Using the inequalities (\ref{z_ratio}) and (\ref{drei}) we obtain

\begin{coro}
\[\lim_{n\to\infty}\frac{\operatorname{mean}(\Delta_n\cdot\gamma)}{\max(\Delta_n\cdot\gamma)}=\frac d{d+1}.\]
\end{coro}

\textbf{Proof of \ref{thm_13}}\begin{enumerate}\item[b)] For the sake of completeness we include a proof of this well known fact.
\begin{equation}\label{Hilb}h_j=\#\{m\in\mathds N^d|m_1+\ldots+m_d=j\}={d-1+j\choose j}\end{equation}
if $\gamma_1=\ldots=\gamma_d=1$. The combinatorial formula
\[\sum_{j=0}^n{d-1+j\choose j}={d+n\choose n}\]
easily yields
\begin{equation}\label{fuenf}n\sum_{j=0}^n{d-1+j\choose j}=\frac{d+1}d\sum_{j=0}^nj{d-1+j\choose j}\end{equation}
and $h(n,\gamma)=\frac d{d+1}$ by (\ref{Hilb}) and (\ref{fuenf}).

\item[a)] Equivalently, we shall prove that for
\[q(n,\gamma):=1-h(n,\gamma)=\frac{\sum_{j=0}^n(n-j)h_j}{n\sum_{j=0}^nh_j},\]
\begin{equation}\label{q_ratio}\lim_{n\to\infty}q(n,\gamma)=\frac1{d+1}.\end{equation}

Proof of (\ref{q_ratio}): Set
\[D:=\{x\in\mathds R_{\geq0}^d|x\cdot\gamma<1\}, D_n:=(n+1)D=\{x\in\mathds R_{\geq0}^d|x\cdot\gamma<n+1\}.\]
Then $(H_0~\dot{\cup}~H_1~\dot{\cup}~\ldots~\dot{\cup}~H_n)\cap\mathds N^d=\Delta_n=D_n\cap\mathds N^d$, hence
\[\sum_{j=0}^nh_j=\#\Delta_n=\#(D\cap\frac1{n+1}\mathds N^d).\]
Analogously, let
\begin{align*}\widehat D&:=\{(x_0,x)\in\mathds R_{\geq0}\times\mathds R_{\geq0}^d|x_0+x\cdot\gamma<1\},\\ \widehat{D_n}&:=(n+1)\widehat D=\{(x_0,x)\in\mathds R_{\geq0}^{d+1}|x_0+x\cdot\gamma<n+1\}.\end{align*}
For the the closures $\overline D$ and $\overline{\widehat D}$ of $D$ and $\widehat D$ we have:

$\overline{\widehat D}$ is a cone over $\overline D$ of height one, hence
\begin{equation}\label{Vol_Cone}\operatorname{Vol}_{d+1}\overline{\widehat D}=\frac1{d+1}\operatorname{Vol}_d\overline D.\end{equation}
Set $\widehat{\Delta_n}:=\widehat{D_n}\cap\mathds N^{d+1}$, hence $\#\widehat{\Delta_n}=\#(\widehat D\cap\frac1{n+1}\mathds N^{d+1})$. By the definition of the Riemann integral
\begin{equation}\label{Vol_1}\lim_{n\to\infty}\frac{\#\Delta_n}{(n+1)^d}=\operatorname{Vol}_d\overline D\end{equation}
and
\begin{equation}\label{Vol_2}\lim_{n\to\infty}\frac{\#\widehat{\Delta_n}}{(n+1)^{d+1}}=\operatorname{Vol}_{d+1}\overline{\widehat D}.\end{equation}
As is easily seen, for $m\in H_j\cap\mathds N^d$, $(y_0,m)\in\widehat{\Delta_n}$ if and only if $0\leq y_0\leq n-j$. Hence above each point $m\in H_j\cap\mathds N^d$ there are exactly $n-j+1$ points of $\widehat{\Delta_n}$. Summing up over all $m\in\Delta_n=\mathop{\dot{\bigcup}}_{j=0}^nH_j\cap\mathds N^d$ we get

\begin{lemma}\begin{equation}\label{Lemma}\#\widehat{\Delta_n}=\sum_{j=0}^n(n-j+1)h_j=\sum_{j=0}^n(n-j)h_j+\#\Delta_n.\end{equation}\end{lemma}

\vspace{.5cm}

Finally, the equations (\ref{Vol_Cone}) up to (\ref{Lemma}) together yield:
\[q(n,\gamma)=\frac{\#\widehat{\Delta_n}-\#\Delta_n}{n\cdot\#\Delta_n}=\frac{\frac{\#\widehat{\Delta_n}}{(n+1)^{d+1}}}{\frac{\#\Delta_n}{(n+1)^d}}\frac{n+1}n-\frac1n\]
converges to $\frac1{d+1}$ for $n\to\infty$.

\end{enumerate}

\hfillplus$\square$

\section{The Hilbert function of positively graded algebras}

Let $\gamma\in\mathds N_{\geq1}^d$ as above. Then for $m=(m_1,\ldots,m_d)\in\mathds N^d$, $m\cdot\gamma$ is the degree of the monomial $X^m=X_1^{m_1}\cdot\ldots\cdot X_d^{m_d}$, hence, with $h_j$ as in section \ref{section_1},
\[f(z):=\sum_{j\geq0}h_jz^j\]
is the Hilbert series of the polynomial ring $\mathds C[X_1,\ldots,X_d]$ w.\,r. to the grading induced by $\deg X_i:=\gamma_i$, $i=1,\ldots,d$. It is well known, that
\[f(z)=\frac1{\prod_{i=1}^d(1-z^{\gamma_i})}.\]
From the results of section 1 we get

\begin{proposition}

\begin{equation}\label{w_ratio}d\geq \frac{\sum_{j=0}^njh_j}{\sum_{j=0}^n(n-j)h_j}=\frac{n\sum_{j=0}^nh_j}{\sum_{j=0}^n(n-j)h_j}-1\buildrel{n\to\infty}\over\longrightarrow d\end{equation}
where $d$ is the pole order at $z=1$ of the rational function $f(z)$.\end{proposition}

Further $f'(z)=\sum_{j\geq1}jh_jz^{j-1}$. Hence (\ref{w_ratio}) may be considered as an analog to the residue formula
\[\frac1{2\pi i}\int_\Gamma\frac{f'(z)}{f(z)}dz=-d\]
where $\Gamma$ is any small circle around the point $z=1$.

A simple but not too simple example of an explicite formula for such a Hilbert function can be found in \cite[4.4.2 Example]{S}: If $d=3$ and $\gamma=(1,2,3)$ then
\[h_n=\lfloor\frac{n^2}{12}+\frac n2\rfloor+1.\]
Using Zhai's lemma cited above and a theorem of Macaulay (see \cite[Theorem 6.1.4]{HH}) it is shown in \cite[3.5 Remark]{HRW2}:

\begin{proposition} \label{Prop_22}Let $I\subseteq\mathds C[X_1,\ldots,X_d]$ be an Artinian ideal, homogenous w.\,r. to the grading induced by $\deg X_i=\gamma_i$, $i=1,\ldots,d$. Then the Hilbert function $(h_n)_{n\in\mathds N}$ of $\mathds C[X_1,\ldots,X_d]/I$ satisfies
\[\sum_{j=0}^mjh_j/\sum_{j=0}^mh_j\leq\frac d{d+1}\cdot m\]
if $h_m\neq0$ and $h_n=0$ for $n>m$.\end{proposition}

\section{Comparison with a question of Wilf}

\label{section_Comparison_Wilf}Let $S$ be a numerical semigroup and $g_0<g_1<\ldots<g_d$ its minimal set of generators and $n_0:=\lceil\frac c{g_0}\rceil$, i.\,e. $0\leq\rho:=n_0g_0-c<g_0$.

Further let $A:=\{s\in S|s-g_0\notin S\}$ the \textbf{Ap\'{e}ry set} of $S$.

Endow $\mathds N^d$ with the \textbf{(purely) lexicographic order} LEX, i.\,e. $a<_{\text{LEX}}b$ if the leftmost nonzero component of $a-b$  is negative. For $a\in A$, let $\tilde a$ be the LEX-minimal element $x\in\mathds N^d$ with $x\cdot(g_1,\ldots,g_d)=a$ and set
\[\tilde A:=\{\tilde a|a\in A\}.\]
In \cite{W}, Wilf raised the following question: Let $\Omega$ be the number of positive integers not contained in $S$ and $c-1$ the largest such element. Is it true, that the fraction $\frac\Omega c$ of omitted numbers is at most $\frac d{d+1}$?

For $d\geq3$ the answer is still unknown. Suppose that the multiplicity $g_0$ of $S$ divides $c$, i.\,e. $c=n_0g_0$.

In this special situation we will see how the problem of Wilf is connected with the considerations on the distribution of weights in section 1. Choose $\gamma_i:=\frac{g_i}{g_0}$, $i=1,\ldots,d$. Then the strips $H_i$ are given by
\[H_i=\{x\in\mathds R^d|ig_0\leq x_1g_1+\ldots+x_dg_d<(i+1)g_0\}, i\in\mathds Z.\]
$\tilde A$, as well as $\Delta_n$ from section 1, is the complement of an $\mathds N^d$-ideal. We define $h_i(S):=\#(H_i\cap\tilde A)$, $i\in\mathds N$. Then, analogously to $\Delta_n$ and $h_i$ in section 1, $h_{n_0}(S)\neq0$, $h_i(S)=0$ for $i\geq n_0+1$ and $\tilde A=\mathop{\dot{\bigcup}}_{i=0}^{n_0}\tilde A\cap H_i$.

According to \cite[End of section 3.1]{HRW2}, then
\begin{equation}\label{Formel_Hilb}h(S):=\frac{\sum_{j=0}^{n_0}jh_j(S)}{n_0\sum_{j=0}^{n_0}h_j(S)}=\frac\Omega c.\end{equation}
Hence, the question of Wilf is, if $h(S)\leq\frac d{d+1}$.

Analogously we may ask:

If $\gamma=\left(\frac{g_1}{g_0},\ldots,\frac{g_d}{g_0}\right)$, do we have (with the notation of section 1)

(i) $h(n,\gamma)\leq\frac d{d+1}$ for $n\geq n_0$?

or

(ii) $h(n_0,\gamma)\leq\frac d{d+1}$?

or at least

(iii) $h(n,\gamma)\leq\frac d{d+1}$ for $n\gg0$?

\section{$A$-graded numerical semigroups}

Let $S$ be a numerical semigroup with minimal set of generators $g_0<g_1<\ldots,<g_d$ and Ap\'{e}ry set $A$. For $i\in\mathds N$ let
\[A_i:=\left\{a\in A\left|\left\lfloor\frac {a+\rho}{g_0}\right\rfloor=i\right.\right\}.\]
We call $S$ $\mathbf A$\textbf{-graded} if $(A_i+A_j)\cap A\subseteq A_{i+j}$, $i,j\in\mathds N$.

\begin{proposition} If $S$ is $A$-graded, then Wilf's inequality $\frac\Omega c\leq \frac d{d+1}$ holds.\end{proposition}

\textbf{Proof }Consider the monomial $\mathds C$-algebra
\[R=R(\tilde A):=\mathds C[X_1,\ldots,X_d]/I, I:=(\underline X^m|m\in\mathds N^d\setminus\tilde A),\]
where $\underline X^m:=X_1^{m_1}\cdot\ldots\cdot X_d^{m_d}$ for $m=(m_1,\ldots,m_d)\in\mathds N^d$.

According to Zhai (\cite{Z}), $\mathds N^d\setminus\tilde A$ is a $\mathds N^d$-ideal, hence $(\underline x^m|m\in\tilde A)$, $\underline x:=\underline X\operatorname{mod}I$, is a $\mathds C$-basis of $R$. Set
\[R_i:=\bigoplus_{a\in A_i}\mathds C\cdot\underline x^{\tilde a}\text{, hence }R=\bigoplus_{i=0}^{n_0}R_i.\]
Since $S$ is $A$-graded, the decomposition $R=\bigoplus_{i=0}^{n_0}R_i$ induces a $\mathds Z$-grading on $R$, that is $R_iR_j\subseteq R_{i+j}$ for $i,j=0,\ldots,n_0$. Further by Zhai \cite[Lemma 1]{Z},
\[\sum_{j=0}^{n_0}jh_j=\sum_{a\in A}\lfloor\frac{a+\rho}{g_0}\rfloor=\Omega+\rho.\]

Applying \ref{Prop_22} Proposition to the positively graded algebra $R=\bigoplus_{i=0}^{n_0}R_i$ yields
\[0\geq(d+1)\sum jh_j-dn_0\sum h_j=(d+1)(\Omega+\rho)-d(c+\rho)\geq(d+1)\Omega-dc.\]
\hfillplus$\square$

\begin{exampleso}

\begin{enumerate}

\item[a)] The semigroup $S=<n^2,n^2+1,n^2+n,n^2+n+1>$, $n\geq3$ is standard $A$-graded, $c=n_0g_0$, $n_0=n-1$, type $t=2n-1$. For different proofs that $S$ satisfies Wilf's condition see also \cite[Thm. 7.1, Cor. 7.2]{E}, \cite[Cor. 2.2]{HRW1} and \cite[2.9 Remark]{HRW2}.

\item[b)] $S=<p,2p+1,2p+3,3p+4>$, $p\geq9$ is $A$-graded with $\deg x_1=\deg x_2=2, \deg x_3=3$, $c=\left\lfloor\frac{2p}3\right\rfloor\cdot p$ and type $t=5$.

\end{enumerate}

\end{exampleso}

\end{document}